\documentclass[12pt,amssymb]{article}%
\usepackage{latexsym}
\usepackage{amsthm}
\usepackage{amsmath}
\usepackage{graphicx}
\usepackage{amsfonts}
\usepackage{amssymb}%
\newtheorem{theorem}{Theorem}
\newtheorem{lemma}{Lemma}
\newtheorem{prop}{Proposition}
\newtheorem{corollary}{Corollary}

\begin{document}

\title{A Dirichlet problem on the half-line for nonlinear equations with indefinite weight}
\author{Zuzana Do\v{s}l\'{a}, Mauro Marini and Serena Matucci}
\maketitle

\begin{center}

\end{center}

\noindent\textbf{Abstract. }{\small We study the existence of positive
solutions on the half-line }$[0,\infty)$ {\small for the nonlinear second
order differential equation }%
\[
\bigl(a(t)x^{\prime}\bigr)^{\prime}+b(t)F(x)=0,\quad t\geq0,
\]
{\small satisfying Dirichlet type conditions, say }$x(0)=0${\small , }%
$\lim_{t\rightarrow\infty}x(t)=0${\small . The function }$b$ {\small is
allowed to change sign and the nonlinearity }$F$ {\small is assumed to be
asymptotically linear in a neighborhood of zero and infinity. Our results
cover also the cases in which }$b$ {\small is a periodic function for large
}$t$ {\small or it is unbounded from below.}

\bigskip

\noindent\textit{Keywords. }{\small Second order nonlinear differential
equation, boundary value problem on the half line, Dirichlet conditions,
globally positive solution, disconjugacy, principal solution.}

{\small \medskip\noindent}

\noindent{\small \textit{MSC 2010:} Primary 34B40, Secondary 34B18.}

\section{Introduction}

Consider the boundary value problem (BVP) on the half-line $[0,\infty)$
\begin{equation}
\bigl(a(t)x^{\prime}\bigr)^{\prime}+b(t)F(x)=0, \label{1}%
\end{equation}%
\begin{equation}
x(0)=0,\quad x(t)>0\text{ \ on }(0,\infty),\text{ }\lim_{t\rightarrow\infty
}x(t)=0, \label{2}%
\end{equation}

where we assume the following:

$\mathrm{{(i)}}$ The function $a$ is continuous on $[0,\infty)$, $a(t)>0,$
and
\begin{equation}
\int_{0}^{\infty}\frac{1}{a(t)}\,dt<\infty. \label{a}%
\end{equation}

$\mathrm{{(ii)}}$ The function $b$ is continuous on $[0,\infty),$ nonnegative
and not identically zero on $[0,1]$, and is allowed to change sign for $t>1$.
Moreover, $b$ is bounded from above, that is, there exists a positive constant
$B$ such that%
\begin{equation}
b(t)\leq B\text{ \ on }[1,\infty). \label{BB}%
\end{equation}

$\mathrm{{(iii)}}$ The function $F$ is continuous on $\mathbb{R}$, $F(u)u>0$
for $u\neq0$, $F$ is differentiable on $[0,\infty)$ with bounded nonnegative
derivative:
\begin{equation}
\ 0\leq\frac{dF(u)}{du}\leq1\ \ \text{ for }u\geq0, \label{DFB}%
\end{equation}
and satisfies
\begin{equation}
\lim_{u\rightarrow0^{+}}\frac{F(u)}{u}=k_{0},\text{ \ \ \ \ }\lim
_{u\rightarrow\infty}\frac{F(u)}{u}=k_{\infty}, \label{FF}%
\end{equation}
where%
\[
0\leq k_{0}\neq k_{\infty}.
\]
Observe that (\ref{DFB}) implies that $k_{0},k_{\infty}\leq1$.

The BVP (\ref{1})-(\ref{2}) is a Dirichlet-type BVP on a unbounded domain.
Recently, there has been a growing interest in studying infinite interval
problems associated to second order nonlinear differential equations, under
various points of view. For a wide bibliography, we refer the reader to
\cite{AO,AGG,kig} and the references therein. When the weight $b$ is of fixed
sign or it is sign-indefinite, we refer to \cite{CDG,MB,WW} or
\cite{CAA,DSA,Trieste,Praga}, respectively. The BVP (\ref{1})-(\ref{2}) arises
in the investigation of positive radial solutions for elliptic equations, when
the nonlinearity is asymptotically linear, see, e.g.,\textbf{ }\cite{BCS}%
\textbf{.}

Our main aim is to continue this study when the function $b$ is allowed to
change its sign and the nonlinearity $F$ can be, roughly speaking, close to a
linear function. The investigated problem can be viewed as an extension to the
half-line of recent results on nonlinear BVPs with a sign-indefinite weight on
a compact interval, see, e.g., \cite{BZ0, BZ1}, and reference therein for a
brief survey on this topic. \vskip2mm Denote by $|\cdot|_{L}$ the norm in
$L^{1}[0,1]$ and set
\begin{equation}
A(t)=\int_{0}^{t}\frac{1}{a(s)}ds. \label{AA}%
\end{equation}
Our main result is the following, in which the disconjugacy of a suitable
auxiliary differential linear equation plays a key role, see Section 3 below.

\begin{theorem}
\label{TMain} Assume that the linear differential equation%
\begin{equation}
v^{\prime\prime}+\frac{B}{a(t)}v=0 \label{DKa}%
\end{equation}
is disconjugate on $[1,\infty),$ where the constant $B$ is defined in
(\ref{BB}).

If there exist $t_{1},t_{2}\in(0,1)$, $t_{1}<t_{2}$ such that $\int_{t_{1}%
}^{t_{2}}b(t)\,dt>0$, and
\begin{align}
&  0\leq\min\{k_{0},k_{\infty}\}A(1)\text{ }|b|_{L}<1,\label{A1}\\[1mm]
&  \max\{k_{0},k_{\infty}\}\int_{t_{1}}^{t_{2}}b(t)\,dt>\frac{A(1)}%
{A(t_{1})(A(1)-A(t_{2}))}, \label{A2}%
\end{align}
then the BVP (\ref{1})-(\ref{2}) has a solution.

Moreover, the solution $x$ has a local maximum in the interval $(0,1]$, is
decreasing in $[1,\infty)$ and satisfies%
\begin{equation}
\int_{1}^{\infty}\frac{1}{a(t)x^{2}(t)}dt=\infty. \label{int}%
\end{equation}

\end{theorem}

\vskip4mm

Theorem \ref{TMain} covers also the cases in which the weight $b$ is a
periodic function for large $t$ or it is unbounded from below.\vskip4mm

Our approach is based on a shooting method and a continuity result. More
precisely, Theorem \ref{TMain} is proved by considering two auxiliary BVPs,
the first one on the compact interval $[0,1]$, where $b$ is nonnegative, and
the second one on the half-line $[1,\infty)$, where $b$ is allowed to change
its sign. The problem of the existence of solutions for (\ref{1}), emanating
from zero, positive in the interval $(0,1)$, and satisfying additional
assumptions at $t=1,$ is considered in Section 2 and is solved by using some
results from \cite{lw98}, with minor changes. The BVP on $[1,\infty)$ is
examined in Section 4. It deals with positive decreasing solutions on
$[1,\infty)$ for (\ref{1}) which tend to zero as $t\rightarrow\infty$. This
second problem is solved by using a fixed point theorem for operators defined
in a Fr\'{e}chet space by a Schauder's linearization device, see \cite[Theorem
1.3]{Furi}. This method does not require the explicit form of the fixed point
operator, but only some \textit{a-priori} bounds. These estimations are
obtained using some properties of principal solutions of disconjugate second
order linear equations, see \cite[Chapter 11]{H}. Finally, roughly speaking,
the solvability of (\ref{1})-(\ref{2}) is obtained by using a shooting method
on $[0,1]$ and, by means some continuity arguments, pasting a solution of
(\ref{1}) on $[0,1]$ with a solution of the BVP on $[1,\infty)$. This last
argument can be viewed as a generalization to non compact intervals of some
ideas in \cite{ghz1}.

Notice that our approach allows us to obtain also an estimation of the decay
to zero of solutions of (\ref{1})-(\ref{2}). Some examples complete the paper.

\section{Two auxiliary BVPs on $[0,1]$}

In this section, we recall some results about the existence of solutions of
(\ref{1}) on $[0,1]$, which belong either to $\Delta_{1}$ or $\Delta_{2}$,
where
\begin{align*}
\Delta_{1}  &  =\left\{  u\in C[0,1]:u(0)=u(1)=0,u(t)>0\text{ on
}(0,1)\right\} \\
\Delta_{2}  &  =\left\{  u\in C[0,1]:u(0)=u^{\prime}(1)=0,u(t)>0\text{ on
}(0,1]\right\}  .
\end{align*}
These results can be obtained from \cite{lw98}, with minor changes.\vskip4mm

BVPs on a compact interval, associated to equations of the form
\begin{equation}
z^{\prime\prime}+g(t)F(z)=0, \label{ea1}%
\end{equation}
where $g$ is a continue nonnegative function on $[0,1]$, have been widely
investigated in the literature, under many different points of view. We refer
to \cite[Introduction]{Calamai} and references therein for a brief survey.

In particular, the existence of solutions of (\ref{ea1}), which satisfy either
$z\in\Delta_{1}$ or $z\in\Delta_{2}$, has been considered in \cite{ew94},
where the key conditions on the nonlinearity are either that $F$ is
superlinear, that is, $k_{0}=0,k_{\infty}=\infty$ or $F$ is sublinear, that
is, $k_{0}=\infty$, $k_{\infty}=0$. When the nonlinearity $F$ is not
necessarily superlinear nor sublinear, these results have been extended in
several ways in \cite{lw98}.

Using \cite[Corollaries 3.1 and 3.5]{lw98} and the continuity of $g,$ we
obtain the following result.

\begin{lemma}
\label{Cor W 1} Assume that there exist $t_{1},t_{2}\in(0,1)$, $t_{1}<t_{2}$,
such that \mbox{$\int_{t_{1}}^{t_{2}}g(t)\,dt>0$} and
\begin{equation}
0\leq\min\{k_{0},k_{\infty}\}\,|g|_{L}<1,\quad\max\{k_{0},k_{\infty}%
\}\int_{t_{1}}^{t_{2}}g(t)\,dt>\frac{1}{t_{1}(1-t_{2})}\, . \label{cor}%
\end{equation}
Then (\ref{ea1}) has both solutions $z_{1}\in\Delta_{1}$ and $z_{2}\in
\Delta_{2}$.
\end{lemma}

\noindent\textit{Proof. } In virtue of the continuity of $g,$ every
nonnegative solution $z$ of (\ref{ea1}), $z\not \equiv 0$, satisfies $z(t)>0$
on $(0,1)$, since $z^{\prime}$ is nonincreasing. Hence, the assertion follows
from \cite[Corollaries 3.1 and 3.5]{lw98}. \hfill$\Box$ \vskip4mm

When $g$ does not have zeros on $[0,1],$ from Lemma \ref{Cor W 1} we obtain
the following.

\begin{lemma}
\label{Cor W 2} Let $g$ be positive on $[0,1]$. If
\begin{equation}
0\leq\min\{k_{0},k_{\infty}\}|g|_{L}<1,\quad\max\{k_{0},k_{\infty}%
\}\,\min_{t\in\lbrack0,1]}g(t)>27, \label{cor1}%
\end{equation}
then (\ref{ea1}) has both solutions $z_{1}\in\Delta_{1}$ and $z_{2}\in
\Delta_{2}$.
\end{lemma}

\noindent\textit{Proof. } Fixed $t_{1},t_{2}\in(0,1)$, $t_{1}<t_{2}$, we have
\[
\int_{t_{1}}^{t_{2}}g(\tau)\,d\tau\geq(t_{2}-t_{1})\min_{t\in\lbrack
0,1]}g(t).
\]
Thus, the second condition in (\ref{cor}) is satisfied if
\[
\max\{k_{0},k_{\infty}\}\min_{t\in\lbrack0,1]}g(t)\geq\frac{1}{t_{1}%
(1-t_{2})(t_{2}-t_{1})}%
\]
for a suitable choice of $t_{1},t_{2}$. Put $\rho=\varrho(t_{1},t_{2}%
)=t_{1}(1-t_{2})(t_{2}-t_{1})$, it is easily checked that $\varrho$ takes its
maximum $1/27$ on the region $0\leq t_{1}<t_{2}\leq1$ when $t_{1}%
=1/3,t_{2}=2/3$. Therefore, the second inequality in (\ref{cor1}) follows.
\hfill$\Box$ \vskip4mm

Define for $t\in\lbrack0,1]$
\begin{equation}
\tau(t)=\frac{A(t)}{A(1)}, \label{defA}%
\end{equation}
where $A$ is given in (\ref{AA}). Thus, $\tau$ maps the interval $[0,1]$ into
itself. Let $x$ be a solution of (\ref{1}) on $[0,1]$ and put $z(\tau
)=x(t(\tau))$, where $t(\tau)$ is the inverse function of $\tau(t).$ Then, $z$
is a solution on $[0,1]$ of
\begin{equation}
\frac{d^{2}z}{d\tau^{2}}+\tilde{b}(\tau)F(z)=0, \label{e:z}%
\end{equation}
where $\tilde{b}(\tau)=A^{2}(1)a(t(\tau))b(t(\tau))$. Vice versa, if $z$ is a
solution of (\ref{e:z}) on $[0,1]$, then $x(t)=z(\tau(t))$ is a solution of
(\ref{1}) on the same interval. Moreover, it is easy to show that $x$ belongs
to $\Delta_{i}$ if and only if $z\in\Delta_{i},i=1,2.$ Hence, Lemmas
\ref{Cor W 1} and \ref{Cor W 2} read for \eqref{1} as follow.

\begin{prop}
\label{Cor W} Assume that one of the following conditions is satisfied.

\begin{itemize}
\item[\textrm{{(i)}}] There exist $t_{1},t_{2}\in(0,1)$, $t_{1}<t_{2}$ such
that $\int_{t_{1}}^{t_{2}}b(t)\,dt>0$, and
\begin{align*}
&  0\leq\min\{k_{0},k_{\infty}\}\,A(1)\,|b|_{L}<1,\\[1mm]
&  \max\{k_{0},k_{\infty}\}\int_{t_{1}}^{t_{2}}b(t)\,dt>\frac{A(1)}%
{A(t_{1})(A(1)-A(t_{2}))}.
\end{align*}

\item[\textrm{{(ii)}}] $b(t)>0$ on $[0,1]$ and
\begin{align*}
&  0\leq\min\{k_{0},k_{\infty}\}\,A(1)\,|b|_{L}<1,\\[2mm]
&  27<\max\{k_{0},k_{\infty}\}\,A(1)\min_{t\in\lbrack0,1]}b(t).
\end{align*}

\end{itemize}

Then (\ref{1}) has both solutions $x_{1}\in\Delta_{1}$ and $x_{2}\in\Delta
_{2}$.
\end{prop}

\noindent\noindent\textit{Proof.} Since
\[
\int_{0}^{1}\tilde{b}(\tau)\,d\tau=A^{2}(1)\int_{0}^{1}b(t(\tau))a(t(\tau
))\,d\tau=A(1)\int_{0}^{1}b(t)\,dt=A(1)\,|b|_{L}%
\]
and%
\[
\int_{t_{1}}^{t_{2}}b(t)\,dt=\frac{1}{A(1)}\int_{\tau_{1}}^{\tau_{2}}\tilde
{b}(\tau)\,d\tau,
\]
where $\tau_{i}=\tau(t_{i})=A(t_{i})/A(1),i=1,2,$ the assertion follows from
Lemmas \ref{Cor W 1} and \ref{Cor W 2}.

\hfill$\Box$ \vskip4mm

Other sufficient conditions for the existence of solutions of (\ref{1}) in the
sets $\Delta_{1}$ and $\Delta_{2}$, can be obtained in a similar way from
other results in \cite{lw98}.\ 

\section{Principal solutions and disconjugacy}

Consider the linear equation
\begin{equation}
\bigl(a(t)y^{\prime}\bigr)^{\prime}+\beta(t)y=0, \label{bet}%
\end{equation}
where $\beta$ is a continuous function for $t\geq T\geq0$. In our study, an
important role is played by the disconjugacy property and the notion of
principal solutions for (\ref{bet}).

We recall that (\ref{bet}) is said to be \textit{disconjugate} on an interval
$I\subset$ $[T,\infty)$ if any nontrivial solution of (\ref{bet}) has at most
one zero on $I$. We refer to \cite{Cop1,H} and references therein for basic
properties of disconjugacy. In particular, the following results will be
useful in the sequel.

\begin{lemma}
\label{LemDIS} Let $T_{1}\geq T.$ The following statements are equivalent.

\textrm{(i$_{1}$)} Equation (\ref{bet}) is disconjugate on $[T_{1},\infty)$;

\textrm{(i$_{2}$)} Equation (\ref{bet}) is disconjugate on $(T_{1},\infty)$;

\textrm{(i$_{3}$)} Equation (\ref{bet}) has a solution without zeros on
$(T_{1} ,\infty)$.
\end{lemma}

\noindent\textit{Proof. } (i$_{1}$)$\Longleftrightarrow$(i$_{2}$). If
(\ref{bet}) is disconjugate on $[T_{1},\infty),$ then it is disconjugate on
$(T_{1},\infty)$. The vice versa follows from \cite[Theorem 2, Chapt.1]{Cop1},
with minor changes. Finally, (i$_{2}$)$\Longleftrightarrow$(i$_{3}$) follows
from \cite[Corollary 6.1]{H}. \hfill$\Box$ \vskip4mm

The concept of principal solution was introduced in 1936 by W. Leighton and M.
Morse and, later on, analyzed by P. Hartman and A. Wintner, see, e.g.,
\cite[Chapter 11]{H}. If (\ref{bet}) is nonoscillatory, then there exists a
solution $u_{0}$ of (\ref{bet}), which is uniquely determined up to a constant
factor by one of the following conditions (in which $u$ denotes an arbitrary
solution of (\ref{bet}), linearly independent of $u_{0}$):
\begin{gather}
\lim_{t\rightarrow\infty}\frac{u_{0}(t)}{u(t)}=0,\label{pp1}\\
\frac{u_{0}^{\prime}(t)}{u_{0}(t)}<\frac{u^{\prime}(t)}{u(t)}\quad\text{ for
large }t,\nonumber\\
\int^{\infty}\frac{dt}{a(t)u_{0}^{2}(t)}=\infty. \label{pp3}%
\end{gather}
The solution $u_{0}$ is called \textit{principal solution }of (\ref{bet}%
)\textit{ }and any solution $u$ of (\ref{bet}), which is linearly independent
of $u_{0},$ is called a \textit{nonprincipal solution }of (\ref{bet}).
Property (\ref{pp1}) is the simplest and most typical property characterizing
principal solutions, because, roughly speaking, it means that the principal
solution is the smallest one\ in a neighborhood of infinity.\vskip2mm

\noindent\textbf{Remark 1. }If (\ref{bet}) is disconjugate on $[T_{1}%
,\infty),T_{1}\geq T,$ then any principal solution of (\ref{bet}) does not
have zeros on $(T_{1},\infty),$ see \cite[Chapter XI, Exercise 6.6]{H}. Thus,
a necessary condition for positiveness of the principal solution on the open
interval $(T,\infty),$ is the disconjugacy of the equation. Nevertheless,
disconjugacy cannot be sufficient for the positiveness of principal solution
on the close half-line $[T,\infty)$, as the following example shows. \vskip2mm

\noindent\textbf{Example 1.} Consider the equation
\begin{equation}
(a(t)y^{\prime})^{\prime}+y=0,\text{ \ }t\geq0, \label{E1}%
\end{equation}
where $a(1)=1$ and
\[
a(t)=\frac{1+t-2e^{t-1}}{1-t}\text{ if }t\neq1.
\]
Hence, $a$ is a positive continuous function on $[0,\infty)$ and (\ref{a})
holds for $a.$ Using (\ref{pp3}), we get that $y_{0}(t)=te^{-t}$ is the
principal solution of (\ref{E1}). Moreover, in view of Lemma \ref{LemDIS},
equation (\ref{E1}) is disconjugate on $[0,\infty)$. \vskip4mm

Consider now the special case $\beta(t)\equiv M>0$ in (\ref{bet}), i.e. the
equation
\begin{equation}
\bigl(a(t)y^{\prime}\bigr)^{\prime}+My=0. \tag{L}\label{La}%
\end{equation}

In view of Example 1, the disconjugacy of (\ref{La}) on $[T,\infty)$ does not
guarantee the positiveness of principal solution at the initial point $t=T.$
To obtain this additional property, consider the so-called \textit{dual}
\textit{equation to (\ref{La}), }that is the equation%
\begin{equation}
v^{\prime\prime}+\frac{M}{a(t)}v=0, \tag{D}\label{Da}%
\end{equation}
which is obtained from (\ref{La}) by the change of variable
$v(t)=a(t)y^{\prime}(t).$ The dual equation has been often used in the
literature for studying oscillatory properties of second order self-adjoint
linear equations, see, e.g., \cite{MON,INT,Sw}, and, for the half-linear case,
\cite{JDE,Doslybook}.

The following necessary and sufficient condition for the disconjugacy of
(\ref{Da}) holds, see also \cite[page 352]{H}.

\begin{lemma}
\label{Lem2} Equation (\ref{Da}) is disconjugate on $[T,\infty)$ if and only
if (\ref{Da}) has a solution $v_{0}$ such that $v_{0}(t)>0$ on $(T,\infty)$
and $v_{0}^{\prime}(t)>0$ on $[T,\infty).$
\end{lemma}

\noindent\textit{Proof. }Assume that (\ref{Da}) is disconjugate on
$[T,\infty)$. From Lemma \ref{LemDIS} there exists a solution $v_{0}$ of
(\ref{Da}) such that $v_{0}(t)>0$ for $t>T.$ Thus, $v_{0}^{\prime}$ is
decreasing for $t>T.$ We claim that $v_{0}^{\prime}(t)>0$ on the whole
interval $[T,\infty).$ By contradiction, if $v_{0}^{\prime}$ has a zero on
$[T,\infty),$ then there exists $t_{1}>T$ such that $v_{0}^{\prime}(t)\leq
v_{0}^{\prime}(t_{1})<0$ for $t\geq t_{1}.$ Integrating this inequality we get
$v_{0}(t)\leq v_{0}(t_{1})+v_{0}^{\prime}(t_{1})(t-t_{1}),$ which gives a
contradiction with the positiveness of $v_{0}$ when $t$ tends to infinity. The
opposite statement follows again in virtue of Lemma \ref{LemDIS}. \hfill$\Box$ \vskip2mm

From Lemma \ref{Lem2} we obtain the following.

\begin{lemma}
\label{Lem main} If (\ref{Da}) is disconjugate on $[T,\infty),$ then
(\ref{La}) has a principal solution $y_{0}$ such that $y_{0}(t)>0$ on
$[T,\infty)$ and $y_{0}^{\prime}(t)<0$ on $(T,\infty)$.
\end{lemma}

\noindent\textit{Proof. } In view of Lemma \ref{Lem2} and the change of
variable $y(t)=v^{\prime}(t),$ equation (\ref{La}) has a solution $y_{0}$
which satisfies $y_{0}(t)>0$ on $[T,\infty)$ and $y_{0}^{\prime}(t)<0$ on
$(T,\infty)$. \ Hence, the disconjugacy of (\ref{La}) follows from Lemma
\ref{LemDIS}. If $y_{0}$ is not principal solution, from \cite[Corollary
6.3]{H} the solution $\overline{y}$ given by%
\[
\overline{y}(t)=y_{0}(t)\int_{t}^{\infty}\frac{ds}{a(s)y_{0}^{2}(s)},
\]
is the desired principal solution of (\ref{La}).\hfill$\Box$ \vskip2mm

\noindent\textbf{Remark 2. }Example 1 shows that the assumption on
disconjugacy of (\ref{Da}) in Lemma \ref{Lem main} cannot by replaced by the
disconjugacy of (\ref{La}). Moreover, observe that the dual equation of
(\ref{E1}) is%
\begin{equation}
v^{\prime\prime}+a^{-1}(t)v=0, \label{ED1}%
\end{equation}
where $a$ is defined in Example 1. It is easy to verify that the function
$v_{0}(t)=2e^{-1}-(1+t)e^{-t}$ is a principal solution of (\ref{ED1}). Since
$v_{0}(1)=0,$ any principal solution of (\ref{ED1}) has a zero at $t=1.$
Consequently, (\ref{ED1}) is not disconjugate on $[0,\infty).$

\section{An auxiliary BVP on $[1,\infty)$}

For any $c>0,$ consider for $t\geq1$ the existence of solutions $x$ of
(\ref{1}) which satisfy the boundary conditions%
\begin{equation}
x(1)=c,\text{ \ }x^{\prime}(1)\leq0,\text{ }x(t)>0\text{ on }[1,\infty),\text{
}\lim_{t\rightarrow\infty}x(t)=0. \label{sec}%
\end{equation}

The solvability of this BVP is based on a general fixed point theorem for
operators defined in a Fr\'{e}chet space, see \cite[Theorem 1.3]{Furi}. In
particular, this result reduces the existence of solutions of a BVP for
differential equations on noncompact intervals to the existence of suitable
\textit{a-priori} bounds and it is mainly useful when the associated fixed
point operator is not known in an explicit form. We recall this result in the
form that will be used.

\begin{theorem}
\label{TeorCFM} Consider the BVP on $[T,\infty),T\geq0,$%
\begin{equation}
(a(t)x^{\prime})^{\prime}+b(t)F(x)=0,\text{ \ \ }x\in S, \label{BVP10}%
\end{equation}
where $S$ is a nonempty subset of the Fr\'{e}chet space $C[T,\infty).$ Let $G$
be a continuous function on $\mathbb{R}^{2},$ such that $F(d)=G(d,d)$ for any
$d\in\mathbb{R}$ and assume that there exist a nonempty, closed, convex and
bounded subset $\Omega\subset C[T,\infty)$ such that for any $u\in\Omega$ the
BVP on $[T,\infty)$
\[
(a(t)x^{\prime})^{\prime}+b(t)G(u(t),x(t))=0,\text{ \ \ }x\in S
\]
admits a unique solution $x_{u}.$ Let $\Psi$ be the operator $\Omega
\rightarrow C[T,\infty),$ such that $\Psi(u)=x_{u}$. Assume

(i$_{1}$) $\Psi(\Omega)\subset\Omega;$

(i$_{2}$) if $\left\{  u_{n}\right\}  \subset\Omega$ is a sequence converging
in $\Omega$ and $\Psi(u_{n})\rightarrow x$, then $x\in S$.

Then $\Psi$ has a fixed point in $\Omega,$ which is a solution of the BVP
(\ref{BVP10}).
\end{theorem}

Let $\widetilde{F}$ be the function%
\begin{equation}
\widetilde{F}(v)=\frac{F(v)}{v}\text{ if }v>0,\text{ \ }\widetilde{F}%
(0)=k_{0}, \label{Fv}%
\end{equation}
where $k_{0}$ is defined in (\ref{FF}) and set $b_{+}(t)=\max\left\{
b(t),0\right\}  ,$ $b_{-}(t)=-\min\left\{  b(t),0\right\}  .$ Thus
$b(t)=b_{+}(t)-b_{-}(t).$ The following holds.

\begin{theorem}
\label{Th second} Assume that equation (\ref{DKa}) is disconjugate on
$[1,\infty)$. Then, for any $c>0,$ equation (\ref{1}) has a unique globally
positive decreasing solution $x$ on $[1,\infty)$ satisfying (\ref{sec}) and
(\ref{int}).
\end{theorem}

\noindent\textit{Proof. } Fixed $c>0,$ consider the equations%
\begin{gather}
\bigl(a(t)y^{\prime}\bigr)^{\prime}+By=0,\label{Ka}\\
\bigl(a(t)w^{\prime}\bigr)^{\prime}-b_{-}(t)w=0. \label{Kb}%
\end{gather}
From Lemma \ref{Lem main}, equation (\ref{Ka}) is disconjugate on $[1,\infty)$
and has a principal solution $y_{0}$ such that $y_{0}(1)=c,$ $y_{0}(t)>0$ on
$[1,\infty),$ $y_{0}^{\prime}(t)<0$ on $(1,\infty).$ Moreover, from
\cite[Theorem 1]{INT} we obtain $\lim_{t\rightarrow\infty}y_{0}(t)=0.$

Since $-b_{-}(t)\leq0,$ equation (\ref{Ka}) is a Sturm majorant for
(\ref{Kb}). Thus (\ref{Kb}) has a positive principal solution $w_{0}$ such
that $w_{0}(1)=c$, $w_{0}^{\prime}(t)\leq0$ for $t\geq1,$ see, e.g.,
\cite[Corollary 6.4]{H}. Using the comparison result for the principal
solutions, see e.g. \cite[Corollary 6.5]{H}, we get on $(1,\infty)$%
\[
\frac{w_{0}^{\prime}(t)}{w_{0}(t)}\leq\frac{y_{0}^{\prime}(t)}{y_{0}(t)}%
\]
and so $0<w_{0}(t)\leq y_{0}(t)$ for $t\geq1$.

Let $\Omega$ and $S$ be the subsets of the Fr\'{e}chet space $C[1,\infty)$
given by
\begin{align*}
\Omega &  =\left\{  u\in C[1,\infty),\frac{1}{2}w_{0}(t)\leq u(t)\leq
y_{0}(t)\right\}  ,\\
S  &  =\left\{  x\in C[1,\infty),x(1)=c,\text{ }x(t)>0,\text{\ }\int
_{1}^{\infty}\frac{1}{a(t)x^{2}(t)}dt=\infty\right\}  ,
\end{align*}
respectively.

For any $u\in\Omega$ consider the linear equation%
\begin{equation}
(a(t)x^{\prime})^{\prime}+b(t)\widetilde{F}(u(t))x(t)=0, \label{Fur}%
\end{equation}
where $\widetilde{F}$ is given in (\ref{Fv}). In view of (\ref{DFB}), we have
$\sup_{v\geq0}\widetilde{F}(v)\leq1.$ Hence, (\ref{Ka}) is a majorant for
(\ref{Fur}). Thus, using again the comparison result \cite[Corollary 6.5]{H},
equation (\ref{Fur}) has a unique positive principal solution $x_{u},$ such
that $x_{u}(1)=c$, and for $t>1$%
\[
\frac{x_{u}^{\prime}(t)}{x_{u}(t)}\leq\frac{y_{0}^{\prime}(t)}{y_{0}(t)}.
\]
Hence, taking into account that $y_{0}$ is decreasing to zero as $t$ tends to
infinity, we get%
\begin{equation}%
\begin{array}
[c]{ll}%
0<x_{u}(t)\leq y_{0}(t) & \text{on }[1,\infty),\\
\lim_{t\rightarrow\infty}x_{u}(t)=0, & x_{u}^{\prime}(t)<0\text{ on }%
(1,\infty)
\end{array}
. \label{est}%
\end{equation}
Thus, for any $u\in\Omega,$ equation (\ref{Fur}) has a solution $x_{u}\in S$,
which is unique in view of (\ref{pp3}).

Denote by $\Psi:\Omega\rightarrow C[1,\infty)$ the operator%
\[
\Psi(u)=x_{u}.
\]
Using again the comparison result \cite[Corollary 6.5]{H} for equations
(\ref{Fur}) and (\ref{Kb}), we obtain for any $u\in\Omega$ and $t\geq1$%
\begin{equation}
\frac{w_{0}^{\prime}(t)}{w_{0}(t)}\leq\frac{x_{u}^{\prime}(t)}{x_{u}(t)}.
\label{cont}%
\end{equation}
Then, in view of (\ref{est}) we get for any $u\in\Omega$ and $t\geq1$
\[
w_{0}(t)\leq x_{u}(t)\leq y_{0}(t),
\]
i.e., the operator $\Psi$ maps $\Omega$ into itself$.$

Now, let $\left\{  u_{n}\right\}  \subset\Omega$ be a sequence converging in
$\Omega$ and $x_{u_{n}}=\Psi(u_{n})\rightarrow x$. Clearly $x(1)=c.$ Since
$\overline{\Psi(\Omega)}\subset\Omega$, we get $x(t)>0.$ Moreover, since
$y_{0}$ is a principal solution of (\ref{Ka}), from (\ref{est}) we obtain%
\[
\int_{1}^{\infty}\frac{1}{a(t)x^{2}(t)}dt\geq\int_{1}^{\infty}\frac
{1}{a(t)y_{0}^{2}(t)}dt=\infty.
\]
Thus, $x\in S$ and, by Theorem \ref{TeorCFM}, there exists a fixed point
$\overline{x}$ of $\Psi$ in $\Omega$. Clearly, $\overline{x}$ is a solution of
(\ref{1}) on $[1,\infty)$ and $\overline{x}(1)=c.$ Since $\overline{x}$ is
also a principal solution of (\ref{Fur}) with $u=\overline{x}$, from
(\ref{est}) we get $\overline{x}(t)>0,$ $\overline{x}^{\prime}(t)<0$ for
$t>1,\overline{x}^{\prime}(1)\leq0$ and $\lim_{t\rightarrow\infty}\overline
{x}(t)=0$. Thus $\overline{x}$ is positive decreasing on $(1,\infty)$ and
satisfies (\ref{sec}) and (\ref{int}).

Finally, it remains to verify that (\ref{1}) has a unique solution which
satisfies (\ref{sec}). Let $x,v$ be two positive solutions of (\ref{1})
defined on $[1,\infty)$ and satisfying (\ref{sec}). In view of the first part
of the proof, we can suppose also that
\begin{equation}
\int_{1}^{\infty}\frac{dt}{a(t)x^{2}(t)}=\infty. \label{yy}%
\end{equation}
Denote by $\Phi(u,v)$ the function $(u\geq0,v\geq0)$
\[
\Phi(u,v)=\left\{
\begin{array}
[c]{lll}%
(F(u)-F(v))/(u-v) & \text{if} & u\neq v\\
&  & \\
dF(u)/du & \text{if } & u=v
\end{array}
\right.
\]
and set $z(t)=x(t)-v(t)$. Thus, $z$ is a solution of the equation%
\begin{equation}
\bigl(a(t)z^{\prime}\bigr)^{\prime}+b(t)\overline{\Phi}(t)z=0, \label{zz}%
\end{equation}
where $\overline{\Phi}(t)=\Phi(x(t),v(t))$. In virtue of (\ref{DFB}), we have%
\[
b(t)\overline{\Phi}(t)\leq B.
\]
Since, from Lemma \ref{Lem main}, equation (\ref{Ka}) is disconjugate on
$[1,\infty),$ the equation (\ref{zz}) is disconjugate on $[1,\infty)$ too$.$
Since $z(1)=0,$ the solution $z$ does not have zeros for $t>1$ and so, without
loss of generality, we can suppose $z(t)>0$ for $t>1.$ Because $\lim
_{t\rightarrow\infty}z(t)=0,$ there exists $t_{1}>1$ such that $z^{\prime
}(t_{1})=0.$ Moreover, taking into account that $x$ satisfies (\ref{yy}) and
$z(t)<x(t)$, we get that $z$ is a principal solution of (\ref{zz}). Using
again the comparison result \cite[Corollary 6.5]{H} for equations (\ref{Ka})
and (\ref{zz}), we obtain for $t>1$\hfill%
\begin{equation}
\frac{z^{\prime}(t)}{z(t)}\leq\frac{y_{0}^{\prime}(t)}{y_{0}(t)}, \label{abs}%
\end{equation}
where $y_{0}$ is the positive decreasing principal solution of (\ref{Ka})
defined in the first part of the proof. Thus, the inequality (\ref{abs}) gives
a contradiction at $t=t_{1}$, because
\[
\frac{y_{0}^{\prime}(t_{1})}{y_{0}(t_{1})}<0.\
\]
\hfill$\Box$ \vskip4mm

We conclude this section with the following continuity result for starting
points of solutions of (\ref{1}) which satisfy (\ref{sec}).

\begin{theorem}
\label{Th cont} Assume that equation (\ref{DKa}) is disconjugate on
$[1,\infty).$ Let $\left\{  c_{n}\right\}  $ be a positive sequence converging
to zero and denote by $x_{n}$ the unique solution of (\ref{1}) which satisfies
(\ref{sec}) with $c_{n}=c.$ Then the sequence $\left\{  x_{n}^{\prime
}(1)\right\}  $ converges to zero.
\end{theorem}

\noindent\textit{Proof. }In virtue of Theorem \ref{Th second}, for any
$c_{n}>0$, equation (\ref{1}) has a unique solution $x_{n}$ which satisfies
(\ref{sec}) with $c_{n}=c.$ Denote by $w_{n}$ the principal solution of
(\ref{Kb}) such that $w_{n}(1)=c_{n}.$ From (\ref{cont}) and (\ref{sec}) we
get%
\begin{equation}
w_{n}^{\prime}(1)\leq x_{n}^{\prime}(1)\leq0. \label{234}%
\end{equation}
Since principal solutions are determined up to a constant factor, we have
\[
w_{n}(t)=\frac{c_{n}}{c_{1}}w_{1}(t).
\]
Hence $w_{n}^{\prime}(1)=c_{n}w_{1}^{\prime}(1)/c_{1}$ and from (\ref{234})
the assertion follows. \hfill$\Box$ \vskip4mm

\section{Proof of the main result}

In this section we prove Theorem \ref{TMain} and we show some its
consequences. To this aim, the following generalization of the well known
Kneser's theorem (see for instance \cite[Section 1.3]{Coppel}), plays a key role.

\begin{prop}
\label{P2} Consider the system
\[
z^{\prime}=F(t,z),\quad(t,z)\in\lbrack T_{1},T_{2}]\times\mathbb{R}^{n}%
\]
where $F$ is continuous and bounded, and let $K_{0}$ be a continuum (i.e., a
compact and connected subset) of $\{(T_{1},w):\,w\in\mathbb{R}^{n}\}$. Let
$\mathcal{Z}(K_{0})$ be the family of all the solutions emanating from $K_{0}%
$. If any solution $z\in\mathcal{Z}(K_{0})$ is defined on the whole interval
$[T_{1},T_{2}]$, then the cross-section $\mathcal{Z}(T_{2};K_{0}%
)=\{z(T_{2}):\,z\in\mathcal{Z}(K_{0})\}$ is a continuum in $\mathbb{R}^{n}$.
\end{prop}

\vskip4mm

\noindent\textit{Proof of Theorem }\ref{TMain}. Let $x_{1}\in\Delta_{1}$ and
$x_{2}\in\Delta_{2}$ be the solutions on $[0,1]$ of (\ref{1}), whose existence
is guaranteed by Proposition \ref{Cor W}, and let $\alpha=\max\{x_{1}^{\prime
}(0),x_{2}^{\prime}(0)\}>0$, $\beta=\min\{x_{1}^{\prime}(0),x_{2}^{\prime
}(0)\}>0$. Put%
\begin{equation}
L=\alpha\,a(0)A(1) \label{LL}%
\end{equation}
and let $\widehat{F}$ be a Lipschitz function on $\mathbb{R}$ such that
\[
\widehat{F}(u)=%
\begin{cases}
0, & u<0\\
F(u), & 0\leq u\leq L\\
F(L), & u>L
\end{cases}
\text{ \ \ .}%
\]
For $\ell\in(0,\alpha],$ consider the Cauchy problem
\begin{equation}%
\begin{cases}
(a(t)x^{\prime})^{\prime}+b(t)\widehat{F}(x)=0, & t\in\lbrack0,1]\\
x(0)=0,\,x^{\prime}(0)=\ell, &
\end{cases}
\label{CPP}%
\end{equation}
and denote by $x_{\ell}$ the unique solution of (\ref{CPP}). Let us show that
$x_{\ell}$ is defined on the whole interval $[0,1]$. For any solution $x$ of
the equation in (\ref{CPP}), the function $a(\cdot)x^{\prime}(\cdot)$ is
nonincreasing, so $a(t)x_{\ell}^{\prime}(t)\leq a(0)x_{\ell}^{\prime
}(0)=a(0)\ell.$ Integrating this inequality, in view of (\ref{LL}) we get for
$t\in\lbrack0,1]$
\[
x_{\ell}(t)\leq a(0)\text{ }\ell\int_{0}^{t}\frac{1}{a(s)}\,ds\leq a(0)\text{
}\ell\text{ }A(1)\leq L.
\]
Assume now that $x_{\ell}(t)>0$ on $(0,t_{1}),0<t_{1}\leq1,$ and $x_{\ell
}(t_{1})=0$. Then, in virtue of the uniqueness of the Cauchy problem
(\ref{CPP}), we obtain $x_{\ell}^{\prime}(t_{1})<0$. If $t_{1}<1$, then
$x_{\ell}(t)<0$ in a right neighborhood of $t_{1}$ and satisfies
$(a(t)x_{\ell}^{\prime})^{\prime}=0$, which gives $x_{\ell}(t)<0$ for every
$t\geq t_{1}$ for which this solution exists. Since $x_{\ell}^{\prime}%
(t_{1})<0,$ by integration we obtain for $t>t_{1}$
\[
x_{\ell}(t)=a(t_{1})x_{\ell}^{\prime}(t_{1})\int_{t_{1}}^{t}\frac{1}%
{a(s)}\,ds>a(t_{1})x_{\ell}^{\prime}(t_{1})A(1),
\]
that is, $x_{\ell}$ is bounded from below. Therefore the solution $x_{\ell}$
of (\ref{CPP}) is defined on the whole interval $[0,1]$.

\smallskip Let $x$ be any solution of (\ref{1}), nonnegative on $[0,1]$ and
satisfying $x(0)=0$, $x^{\prime}(0)=\ell\in(0,\alpha]$. Then $x$ is also a
solution of (\ref{CPP}) for $0\leq t\leq1$, and vice versa. Indeed, reasoning
as above, we obtain $x(t)\leq L$ on $[0,1]$ and therefore $F(x(t))=\widehat
{F}(x(t))$ for all $t\in\lbrack0,1]$.

\smallskip Put $K_{0}=\{(x(1),x^{\prime}(1)):x\text{ is solution of }%
$(\ref{CPP})$\text{ with }\ell\in\lbrack\beta,\alpha]\}$. Since any solution
of (\ref{CPP}) is defined on the whole $[0,1]$, by Proposition \ref{P2} the
set $K_{0}$ is a continuum in $\mathbb{R}^{2}$, containing the points
$(0,x_{1}^{\prime}(1))$, $(x_{2}(1),0)$, with $x_{1}^{\prime}(1)<0$,
$x_{2}(1)>0$. Further, $K_{0}$ does not contain any point $(0,c)$ with
$c\geq0$. Therefore a continuum $K_{1}\subseteq K_{0}$ exists, $K_{1}%
\subseteq\overline{\pi}=\{(u,v):\,u\geq0,v\leq0\}$, $(0,0)\notin K_{1}$, and
there exist two points $P,Q\in K_{1},$ $P=(p,0),$ $Q=(0,-q),$ $p>0,q>0$.

In order to complete the proof, we use a similar argument to the one given in
\cite[Theorem 1.1]{Trieste}, with minor changes. \ Consider equation (\ref{1})
for $t\geq1$. By Theorem \ref{Th second}, for every $c>0$, (\ref{1}) has a
unique positive decreasing solution $x$ satisfying (\ref{sec}) and
(\ref{int}). Then, the set $S_{1}$ of the initial data of the solutions of
(\ref{1}) on $[1,\infty)$ satisfying (\ref{int}) and (\ref{sec}) is connected,
$S_{1}\subset\bar{\pi}$, and its projection on the first component is the
half-line $(0,\infty)$. Further, from Theorem (\ref{Th cont}), $(0,0)\in
\bar{S}_{1}$. Therefore we have
\[
K_{1}\cap S_{1}\neq\emptyset.
\]
Let us show that to each point $(c,d)\in K_{1}\cap S_{1}$ corresponds to a
solution of the BVP (\ref{1})-(\ref{2}). Let $(c,d)\in K_{1}\cap S_{1}$. Then
$c>0,d\leq0$. Since $(c,d)\in K_{1}$, there exists a solution $u$ of
(\ref{CPP}), for a suitable $\ell\in\lbrack\beta,\alpha]$, such that
$u(1)=c>0$ and $u^{\prime}(1)=d\leq0$. Since $u(1)>0$ we have $u(t)>0$ on
$(0,1]$. Therefore $u$ is also a solution of (\ref{1}) in $[0,1]$, with
$u(0)=0$, $u(t)>0$ for $t\in(0,1]$. As $(c,d)\in S_{1}$, a positive decreasing
solution $v$ of (\ref{1}) exists on $[1,\infty)$, which satisfies (\ref{sec})
and $v(1)=c=u(1)$, $v^{\prime}(1)=d=v^{\prime}(1)$. Hence, the function
\[
x(t)=%
\begin{cases}
u(t), & t\in\lbrack0,1],\\
v(t), & t>1.
\end{cases}
\]
is a solution of the BVP (\ref{1})-(\ref{2}) and the proof is complete.\hfill
$\Box$ \vskip4mm\bigskip

From Theorem \ref{TMain} and Proposition \ref{Cor W}, we get the following.

\begin{corollary}
\label{C1} Let assumptions of Proposition \ref{Cor W}-(ii) be satisfied and
equation (\ref{DKa}) is disconjugate on $[1,\infty).$ Then the BVP
(\ref{1})-(\ref{2}) has a solution.
\end{corollary}

\vskip2mm

We close this section with the solvability of our BVP for the perturbed
equation
\begin{equation}
\bigl(a(t)z^{\prime}\bigr)^{\prime}+(b(t)+b_{1}(t))F(z)=0, \label{100}%
\end{equation}
where $b_{1}$ is a continuous function for $t\geq0$ such that $b_{1}%
(t)\equiv0$ on $[0,1]$ and $b_{1}(t)\leq0$ for $t>1.$

\begin{corollary}
\label{C2} If assumptions of Theorem \ref{TMain} are satisfied, then equation
(\ref{100}) has a solution $z$ satisfying boundary conditions (\ref{2}).
\end{corollary}

\section{Examples and concluding remarks}

Theorem \ref{TMain} is illustrated by the following example. \vskip2mm

\noindent\textbf{Example 2.} Consider the equation
\begin{equation}
\bigl(a(t)x^{\prime}\bigr)^{\prime}+b(t)\text{ }F(x)=0, \label{Ex 2.3}%
\end{equation}
where%
\begin{equation}
a(t)=(1+t)^{2},\quad b(t)=\frac{1}{5e}\exp\left(  \frac{16}{1+16t^{4}}\right)
\cos\left(  \frac{\pi t}{2}\right)  \quad\text{for }t\geq0. \label{ab}%
\end{equation}
and $F$ satisfies (\ref{DFB}) and (\ref{FF}) with
\[
k_{0}=\frac{9}{e^{15}},\quad k_{\infty}=1.
\]
Since $b$ is decreasing on $[0,1],$ we get%
\[
\int_{1/3}^{1/2}b(\tau)d\tau\geq\frac{1}{6}b(1/2)=\frac{\sqrt{2}}{60}e^{7}.
\]
For equation (\ref{Ex 2.3}), the function $A$ in (\ref{AA}) becomes%
\[
A(t)=\frac{t}{1+t},
\]
so assumptions (\ref{A1}), (\ref{A2}) are verified for $t_{1}=1/3$ and
$t_{2}=1/2,$ because
\begin{gather*}
A(1)|b|_{L}\geq\frac{1}{2}b(0)=\frac{e^{15}}{10},\\
\frac{A(1)}{A(t_{1})(A(1)-A(t_{2}))}=12<\frac{\sqrt{2}}{60}e^{7}.
\end{gather*}
Finally, for $t\in\lbrack1,\infty)$ we have%
\[
b(t)\leq\frac{1}{5e}<\frac{1}{4}%
\]
and the equation (\ref{DKa}) becomes the Euler equation%
\[
v^{\prime\prime}+\frac{1}{4(1+t)^{2}}v=0,
\]
which is disconjugate on $[1,\infty),$ see, e.g., \cite[Chapter 2.1]{Sw}.
Hence, in view of Theorem \ref{TMain}, equation (\ref{Ex 2.3}) has solutions
$x$ which satisfy the boundary conditions (\ref{2}) and%
\[
\int_{1}^{\infty}\frac{1}{(1+t)^{2}x^{2}(t)}dt=\infty.
\]

\noindent{\textbf{Remark 3.}} Example 2 can be slightly modified for the
nonlinearity
\[
F(u)=\frac{u^{2}}{1+u}%
\]
or the nonlinearity%
\[
F(u)=\frac{u}{1+\sqrt{u}}.
\]
\vskip4mm

\noindent{\textbf{Remark 4.}} Consider the equation%
\begin{equation}
\bigl(a(t)x^{\prime}\bigr)^{\prime}+(b(t)\text{ }+b_{1}(t))F(x)=0,
\label{Ex 2.2}%
\end{equation}
where the functions $a,b$ are given in (\ref{ab}), $b_{1}$ is the function
\[
b_{1}(t)=(e-e^{t})(|\cos t|-\cos t),
\]
and $F$ is as in Example 2. Since $b_{1}(t)\leq0,$ in view of Corollary
\ref{C2}, equation (\ref{Ex 2.2}) has solutions $x$ which satisfy the boundary
conditions (\ref{2}).\vskip4mm

\noindent\textbf{Remark 5.} Theorem \ref{TMain} and Corollaries \ref{C1},
\ref{C2} continue to hold if the assumption (\ref{DFB}) is replaced by the
more general condition%
\[
\exists K>0:0\leq\frac{dF(u)}{du}\leq K\ \ \text{ for }u\geq0
\]
and the disconjugacy of (\ref{DKa}) is substituted by the disconjugacy on
$[1,\infty)$ of the linear equation%
\[
v^{\prime\prime}+\frac{B\,K}{a(t)}v=0.
\]

\noindent\textbf{Remark 6. }The assumption $k_{0}\neq k_{\infty}$ implies that
$F$ cannot be a linear function on $[0,\infty).$ If the linear equation
\begin{equation}
\bigl(a(t)x^{\prime}\bigr)^{\prime}+b(t)x=0 \label{LLL}%
\end{equation}
has a solution $x$ satisfying (\ref{2}), then in virtue of Lemma \ref{LemDIS},
(\ref{LLL}) is disconjugate on $[0,\infty)$. However, $x$ is not necessarily
the principal solution of (\ref{LLL}). The following example illustrates this case.

\vskip1mm

\noindent\textbf{Example 3.} Consider the equation%
\begin{equation}
\bigl(e^{2t}x^{\prime}\bigr)^{\prime}+e^{2t}x=0,\quad t\geq0 \label{Ex L}%
\end{equation}
A standard calculation shows that $x_{0}(t)=e^{-t}$ , $x_{1}(t)=te^{-t}$ are
solutions of (\ref{Ex L}). Obviously, $x_{1}$ satisfies (\ref{2}). Observe
that $x_{1}$ is a nonprincipal solution and $x_{0}$ is the principal solution.

\bigskip

\vskip4mmIn a forthcoming paper we will consider this kind of BVPs for
nonlinear equations for which $k_{0}=k_{\infty}$.

\bigskip

\vskip4mm

\noindent{\small \emph{Authors' addresses}: }

{\small \noindent\emph{Zuzana Do\v{s}l\'{a},} Department of Mathematics and
Statistics, Masaryk University, Kotl\'{a}\v{r}sk\'{a} 2, CZ-61137 Brno, Czech
Republic. E-mail: \texttt{dosla@math.muni.cz}} \vskip2mm

{\small \noindent\emph{Mauro Marini,} Department of Mathematics and
Informatics "Ulisse Dini", University of Florence, I-50139 Florence, Italy.
E-mail: \texttt{mauro.marini@unifi.it}} \vskip2mm

{\small \noindent\emph{Serena Matucci,} Department of Mathematics and
Informatics "Ulisse Dini", University of Florence, I-50139 Florence, Italy.
E-mail: \texttt{serena.matucci@unifi.it}}

\vskip4mm

\noindent

\begin{thebibliography}{99}                                                                                               %


\bibitem {AO}Agarwal R.P., O'Regan D.: Infinite interval problems for
differential, difference and integral equations, Kluwer Academic Publishers,
Dordrecht, 2001.

\bibitem {AGG}Andres J., Gabor G., G\'{o}rniewicz L.: \textit{Boundary value
problems on infinite intervals}, Trans. Amer. Math. Soc. \textbf{351} (1999), 4861--4903.

\bibitem {BCS}Bartolo R., Candela A. M., Salvatore A.: \textit{Perturbed
asymptotically linear problems}, Ann. Mat. Pura Appl. (4) \textbf{193} (2014), 89--101.

\bibitem {BZ0}Boscaggin A., Zanolin F.: \textit{Pairs of positive periodic
solutions of second order nonlinear equations with indefinite weight,} J.
Differential Equations \textbf{252} (2012), 2900--2921.

\bibitem {BZ1}Boscaggin A., Zanolin F.: \textit{Second-order ordinary
differential equations with indefinite weight: the Neumann boundary value
problem}, Ann. Mat. Pura Appl. (2013), DOI 10.1007/s10231-013-0384-0.

\bibitem {Calamai}Calamai A., Infante G.: \textit{Nontrivial solutions of
boundary value problems for second-order functional differential equations},
Ann. Mat. Pura Appl. (2013), DOI 10.1007/s10231-015-0487-x.

\bibitem {CDG}Cort\'{a}zar C., Dolbeault J., Garc\'{\i}a-Huidobro M.,
Man\'{a}sevich R.: \textit{Existence of sign changing solutions for an
equation with a weighted p-Laplace operator}, Nonlinear Anal. \textbf{110}
(2014), 1--22.

\bibitem {MON}Cecchi M., Marini M., Villari Gab.: \textit{On monotonicity
property for a certain class of second order differential equations}, J.
Differential Equations \textbf{82} (1989), 15-27.

\bibitem {INT}Cecchi M., Marini M., Villari Gab.: \textit{Integral criteria
for a classification of solutions of linear differential equations}, J.
Differential Equations \textbf{99} (1992), 381-397.

\bibitem {JDE}Cecchi M., Do\v{s}l\'{a} Z., Marini M.: \textit{Half-linear
equations and characteristic properties of the principal solution}, J.
Differential Equations \textbf{208} (2005), 494-507, Corrigendum J.
Differential Equations \textbf{221} (2006) 272-274.

\bibitem {Furi}Cecchi M., Furi M., Marini M.: \textit{ On continuity and
compactness of some nonlinear operators associated with differential equations
in noncompact intervals,} Nonlinear Anal., T.M.A, \textbf{9} (1985), 171-180.

\bibitem {Coppel}Coppel W.A.: Stability and Asymptotic Behavior of
Differential Equations, D. C. Heath and Co., Boston 1965.

\bibitem {Cop1}Coppel W.A.: Disconjugacy, Lecture Notes Math. \textbf{220},
Springer-Verlag, Berlin, 1971.

\bibitem {MB}Do\v{s}l\'{a} Z., Marini M., Matucci S.: \textit{On some boundary
value problems for second order nonlinear differential equations}, Math.
Bohem. \textbf{137} (2012), 113-122.

\bibitem {CAA}Do\v{s}l\'{a} Z., Marini M., Matucci S.: \textit{A boundary
value problem on a half-line for differential equations with indefinite
weight}, Commun. Appl. Anal. \textbf{15} (2011), 341-352.

\bibitem {DSA}Do\v{s}l\'{a} Z., Marini M., Matucci S.: \textit{Positive
solutions of nonlocal continuous second order BVP's}, Dynamic Syst. Appl.
\textbf{23} (2014), 431-446.

\bibitem {Doslybook}Do\v{s}l\'{y} O., \v{R}eh\'{a}k P.: Half-linear
Differential Equations, North-Holland, Mathematics Studies \textbf{202},
Elsevier Sci. B.V., Amsterdam, 2005.

\bibitem {ew94}Erbe L.H., Wang H.: \textit{On the existence of positive
solutions of ordinary differential equations}, Proc. Amer. Math. Soc.
\textbf{120} (1994), 743-748.

\bibitem {ghz1}Gaudenzi M., Habets P., Zanolin F.: \textit{An example of a
superlinear problem with multiple positive solutions,} Atti Sem. Mat. Fis.
Univ. Modena \textbf{51} (2003), 259--272.

\bibitem {H}Hartman P.: Ordinary Differential Equations, 2 Ed., Birk\"{a}user,
Boston-Basel-Stuttgart 1982.

\bibitem {kig}Kiguradze I.T., Chanturia A.: Asymptotic Properties of Solutions
of Nonautonomous Ordinary Differential Equations, Kluwer Acad. Publ. G.,
Dordrecht, 1993.

\bibitem {lw98}Lan K., Webb J.R.L.: \textit{Positive solutions of semilinear
differential equations with singularities}, J. Differential Equations
\textbf{148} (1998), 407--421.

\bibitem {Trieste}Marini M., Matucci S.: \textit{ A boundary value problem on
the half-line for superlinear differential equations with changing sign
weight}, Rend. Istit. Mat. Univ. Trieste \textbf{44} (2012), 117-132.

\bibitem {Praga}Matucci S.: \textit{A new approach for solving nonlinear BVP's
on the half-line for second order equations and applications,} to appear on
Mathematica Bohemica \textbf{140} (2015). \ \ \ \ \ \ 

\bibitem {Sw}{Swanson C.A.: }Comparison and Oscillation Theory of Linear
Differential Equations, Academic Press, New York, 1968{\textit{\textrm{.}}}

\bibitem {WW}Wei Y., Wong P.J.Y.: \textit{Existence and uniqueness of
solutions for delay boundary value problems with p-Laplacian on infinite
intervals}, Bound. Value Probl. \textbf{2013-141} (2013), 1-13.
\end{thebibliography}
\end{document}